\documentclass[14pt]{article}

\usepackage{multicol}
\usepackage{amsthm}

\usepackage[dvips]{graphicx,color}
\usepackage{amssymb}
\usepackage{amsmath}

\begin{document}

\fontsize{14pt}{16.5pt}\selectfont
~~\\
\vspace{-5mm}

\noindent
\begin{center}

\bf{On a difference in topological nature between Cantor middle-third set and Sierpi\'nski carpet}
\end{center}

\fontsize{12pt}{15pt}\selectfont

Akihiko Kitada$^{1}$,~Tomoyuki Yamamoto$^{1,2}$,~Shousuke Ohmori$^{2}$,~Yoshihiro Yamazaki$^{1,2}$\\

\noindent
$^1$\it{Institute of Condensed-Matter Science, Comprehensive Resaerch Organization, Waseda University,
3-4-1 Okubo, Shinjuku-ku, Tokyo 169-8555, Japan}\\
$^2$\it{Faculty of Science and Engineering, Waseda~University, 3-4-1 Okubo, Shinjuku-ku, Tokyo 169-8555, Japan}\\
~~\\

\rm
\noindent
{\bf{Abstract}}\\
\noindent
Both Cantor middle-third set and Sierpi\'nski carpet are self-similar, perfect, compact metric spaces. In spite of the similarity of the mathematical procedure of construction, there exists between them a fundamental difference in topological nature, and this difference affects the methods of construction of an interesting non-trivial quotient space of them. The totally disconnectedness (or, more generally, zero-dimensional) enables Cantor middle-third set to have a non-trivial quotient space which is self-similar. On the other hand, concerning Sierpi\'nski carpet, because of the connectedness of its structure, no non-trivial quotient space which is self-similar can be constructed by such an elegant procedure as that for Cantor middle-third set. Various topologically significant nature specific to Cantor middle-third set owe mainly to the totally disconnectedness of the set.

\bigskip



\section{Introduction}
\noindent
~~Since the procedure of the block construction or the coarse graining of a space corresponds mathematically to that of the construction of a quotient space (2-D.) of a space, nonlinear science researchers have been interested in quotient spaces[1,2].\\
~~In the present paper, between Cantor middle-third set (CMTS) and Sierpi\'nski carpet (SC) (2-C.) both of which are perfect (2-A.), compact and self-similar (2-B.) a difference in topological nature will be made clear through the investigations of the procedure of construction of an interesting quotient space. In fact, owing to the zero-dimensional (0-dim) (2-A.) or totally disconnectedness of CMTS, some interesting non-trivial quotient space of CMTS such as one with a self-similarity can be constructed in an elegant way. Contrary to this, such an interesting non-trivial quotient space is hardly constructed for SC at least in the same elegant procedure with that for CMTS because of the connectedness of SC. In the following sections, this fact will be confirmed.
 
\section{Definitions and Statements}
~~Some definitions and statements which will be used in the section 3 are listed below.

\begin{description}
\item[2-A.] 
~A topological space $(Z,\tau)$ is said to be zero-dimensional (0-dim) provided that at every point $p\in Z$ and for every $U\in\tau$ containing $p$, there exists $u\in \tau\cap \Im$ ($\Im$ denotes the set of all closed sets of $Z$, that is, $u$ is an open and closed set) such that $p\in u\subset U$. In any compact T$_2$-space, the notion of 0-dim coincides with that of the totally disconnectedness[3,4] and it is easy to check that any 0-dim T$_0$-space is totally disconnected.\\
~~A topological space $(Z,\tau)$ is said to be perfect provided that a singleton $\{p\}$, $p\in Z$ is not $\tau$-open, that is, any non-empty open set of $(Z,\tau)$ has at least two points. Any non-empty open set of a perfect space is perfect as a subspace.\\
~~Both 0-dim and perfect are topological properties, that is, 0-dim and perfect are invariant under the homeomorphism[5].\\
\item[2-B.]
~A topological space $(A,\tau)$ is self-similar provided that i) it is metrizable, that is, $\tau$ is identical with a metric topology $\tau_d$ defined by a metric $d$ on $A$. ii) there exists a set of contractions $\{f_j:(A,\tau_d)\rightarrow (A,\tau_d);j=1,\dots,m~(2\le m<\infty)\}$ such that $\displaystyle\bigcup_{j=1}^mf_j(A)=A $.\\
~~Especially, we call a metric space $(A,\tau_d)$ self-similar provided that the condition ii) is satisfied with respect to this metric $d$.\\
~~A sufficient condition for a topological space to be self-similar is given in the following statement.\\
\end{description}

{\noindent\bf Statement.}~{\it The existence of a self-similar space which is homeomorphic to $(Y,\tau)$ is sufficient for a topological space $(Y,\tau)$ to be self-similar.}\\
~~(Concerning the proof, see Appendix 1).)
\bigskip

\begin{description}
\item[2-C.] 
~CMTS is a subspace of closed interval $[0,1]$ such that CMTS=$\displaystyle\bigcup_{j=1}^2 f_j($CMTS) where the contractions $f_1(x)$ and $f_2(x)$ are $x/3$ and $x/3+2/3,~x\in[0,1]$, respectively. The sum of two contraction coefficients is $2/3<1$.\\~~Any metric space which is 0-dim, perfect and compact is homeomorphic to CMTS[6].\\
~~SC is a subspace of square $[0,1]\times [0,1]$ such that SC=$\displaystyle\bigcup_{j=1}^8 f_j($SC) where\\ $f_1(x_1,x_2)=(x_1/3,x_2/3),~f_2(x_1,x_2)=(x_1/3+1/3,x_2/3),\\
f_3(x_1,x_2)=(x_1/3+2/3,x_2/3),~f_4(x_1,x_2)=(x_1/3,x_2/3+1/3),~f_5(x_1,x_2)=(x_1/3+2/3,x_2/3+1/3),~f_6(x_1,x_2)=(x_1/3,x_2/3+2/3),~f_7(x_1,x_2)=(x_1/3+1/3,x_2/3+2/3),~f_8(x_1,x_2)=(x_1/3+2/3,x_2/3+2/3),\\(x_1,x_2)\in[0,1]\times [0,1]$. The sum of eight contraction coefficients is $8/3\geq 1$. \\
In general, the following statement holds.\\
\end{description}

{\noindent\bf Statement.} ~{\it Let $(Z,\tau_d)$ be a compact metric space and the sytem\\
$\{f_j:(Z,\tau_d)\to(Z,\tau_d),~j=1,\ldots,m\}$ be a set of contractions \\
$d(f_j(z),f_j(z'))\le\alpha_j(\eta)d(z,z')~{\rm for}~d(z,z')<\eta,~0<\alpha_j(\eta)<1,$
$\inf_{\eta>0}\alpha_j(\eta)>0, j= 1, \cdots, m~(2\le m<\infty)$.\\
Then, there exists a non-empty compact set $X\subset Z$ such that $\displaystyle\bigcup_{j=1}^m f_j(X)=X$. This self-similar set $X$ is unique in the set of all non-empty compact sets of $Z$.\\
~~~~If the system of contractions $\{f_j\}$ satisfies two conditions
\begin{itemize}
\item[i)] Each $f_j$ is one to one,
\item[ii)] The set $\displaystyle\bigcup_{j=1}^m \{z\in Z;~f_j(z)=z\}$ is not a singleton,
\end{itemize}
then the self-similar compact set $X$ is perfect.\\
~~~~Furthermore, if the system of contractions $\{f_j\}$ satisfies the additional third condition
\begin{itemize}
\item[iii)] $\displaystyle\sum_{j=1}^m inf_{\eta>0} \alpha_j(\eta)<1$\\
\end{itemize}
the self-similar compact set $X$ is not only perfect but 0-dim.}\\
~~(Concerning the proof, see Appendix 2).)\\
~~\\

Thus, from this statement, SC is perfect and compact, on the other hand, CMTS is perfect, compact and 0-dim.\\

\begin{description}
\item[2-D.]
~A quotient space of a space $X$ is defined by a classification of all points in $X$ through the identificatioin of the different points based on an equivalence relation.\\
~If $f:(X,\tau)\rightarrow (Y,\tau')$ is a continuous, onto, closed map (A map $f$ is said to be a closed map provided that for any closed set $K$ in $X$, $f(K)$ is a closed set in $Y$.), then 
$h:(Y,\tau')\rightarrow (\mathcal{D}_f,\tau(\mathcal{D}_f)),~y\mapsto f^{-1}(y)$ is a 
homeomorphism. Here the set $\mathcal{D}_f$ is $\{f^{-1}(y)\subset X;y\in Y\}$ and the topology $\tau(\mathcal{D}_f)$ is $\{\mathcal{U}\subset \mathcal{D}_f;\displaystyle\bigcup \mathcal{U}=\bigcup _{D\in\mathcal{U}}D\in\tau,~$each$~D\in \mathcal{D}_f$ is $f^{-1}(y)$ for some $y\in Y\}$. Topological space $(\mathcal{D}_f,\tau(\mathcal{D}_f))$ is a quotient space of $(X,\tau)$ based on the equivalence relation $\sim $
\begin{center}
$x\sim x' \overset{\mathrm{def}}{=}
 f(x)=f(x')$
\end{center}
Sometimes, the quotient space $(\mathcal{D}_f,\tau(\mathcal{D}_f))$ is called a decomposition space of $(X,\tau)$ due to $f$.\\
~In general, a quotient space $(\mathcal{D},\tau(\mathcal{D}))$ of $(X,\tau)$ where $\mathcal{D}=\{\{x\};x\in X\}$ ($\{x\}$ denotes a singleton) and $\tau(\mathcal{D})=\{\mathcal{U}\subset \mathcal{D};\displaystyle\bigcup \mathcal{U}\in\tau\}$ is called the trivial quotient space. The trivial quotient space $\mathcal{D}$ of $X$ is homeomorphic to $X$. Therefore, if $X$ is self-similar, from the statement in 2-B., $\mathcal{D}$ must be self-similar.\\
\item[2-E.~Statement.] {\it ~Let $(Z,\tau)$ be a 0-dim, perfect T$_0$(necessarily T$_2$)-space. Then, there exists a 
partition $\{Y, Z-Y\}$ of $Z$ where $Y\in (\tau\cap \Im)-\{\phi \}$ and $Y\subsetneqq Z$. Here the subspace $(Y,\tau_Y)$ is 0-dim and perfect.}\\
~~(Concerning the proof, see Appendix 3).)\\
\bigskip

\item[2-F.][7]
Let $(X,\tau)$ be a compact normal space and $(X_i,\tau_{X_i}),~X_i \in \Im - \{\phi \},~i\in \bf{N}$ be connected subspaces of $(X,\tau)$ such that $X_{i+1}\subset X_i,~i\in \bf{N}$. Then, $E=\displaystyle\bigcap _{i\in \bf{N}}X_i$ is not empty and the subspace $(E,\tau_E)$ of $(X,\tau)$ is compact and connected.
\end{description}
~~\\

\section{Quotient spaces}
\subsection{Quotient space of CMTS}

Since CMTS is a 0-dim, perfect, compact metric space (2-C.), from 2-E., there exists a 0-dim, perfect, compact subspace $(Y,\tau_{d_Y})$ of CMTS such that $Y\subsetneqq$ CMTS. For example, $Y=$CMTS$\cap [0,1/3]$. Here the notation $d_Y$ denotes the restriction of the metric $d$ ($d(x,x')=|x-x'|$) on $Y$. Now, let $f:($CMTS,$\tau_d)\rightarrow (Y,\tau_{d_Y})$ be a not one to one continuous, onto map. The simple one  
\begin{center}~~~
$f(x)=\left\{
\begin{array}{lcr}
x~,~~x\in Y \\
q~,~~x\in $ CMTS $-Y
\end{array}
\right.$
\end{center}
is an example. Here $q$ is an arbitrarily chosen point of $Y$. Since the singleton 
$\{q\}$ is a closed set in a T$_1$-space, it is easy to see that the map $f$ is a 
continuous, closed map from CMTS onto $Y$. From 2-D., $h:(Y,\tau_{d_Y})\rightarrow 
(\mathcal{D}_f,\tau(\mathcal{D}_f)), y\mapsto f^{-1}(y)$ is a homeomorphism. Here, $(\mathcal{D}_f,\tau(\mathcal{D}_f))$ is a quotient space of CMTS due to the above map $f$. Since $(Y,\tau_{d_Y})$ is a 0-dim, perfect, compact metric space, from 2-C., it is homeomorphic to CMTS. Then, the quotient space $\mathcal{D}_f$ is also homeomorphic to CMTS. According to the statement in 2-B., $\mathcal{D}_f$ is self-similar. It must be noted that this $\mathcal{D}_f$ is not the 
trivial quotient set $\{\{x\};x\in $CMTS$\}$ because $f$ is not one to one. In fact, for example, if $Y=$CMTS$\cap [0,1/3]$, the set $\{q\}\cup ($CMTS$\cap [2/3,1])$ is a point of $\mathcal{D}_f$. Namely, there exists a self-similar non-trivial quotient space of CMTS. Furthermore, since $\mathcal{D}_f$ is homeomorphic to CMTS, it is a self-similar, 0-dim, perfect, compact metric space (Concerning the metrizability of $\mathcal{D}_f$, see Appendix 1).). Then, we can obtain a self-similar quotient space of $\mathcal{D}_f$ which is a 0-dim, perfect, compact metric space again. Continuing this process, we attain an infinite sequence $\{\mathcal{D}_i;i=0,1,2,\cdot \cdot \cdot \}$ in which $\mathcal{D}_0=$CMTS and each $\mathcal{D}_i$ is a quotient space of $\mathcal{D}_{i-1}$\\ 

\subsection{Quotient space of SC}

According to the statement in 2-C., SC is perfect and compact. CMTS is also perfect and compact, and furthermore, 0-dim. In this paragraph, using 2-F., we will demonstrate that SC is a connected space. SC is given by the intersection $\displaystyle\bigcap _{i\in \bf{N}}X_i$ of closed sets $X_i=\displaystyle\bigcup_{j{_1}\cdot \cdot \cdot j{_i}}f_{j_1}\circ \cdot \cdot \cdot \circ f_{j_i}([0,1]\times [0,1]),~j_1\in \{1,\dots,8\},\dots,j_i\in \{1,\dots,8\},~i\in \bf{N}$ in $[0,1]\times [0,1]$. Here each $X_i$ has finite number of windows and obeys the relation $X_{i+1}\subset X_i$. Fig.1 shows, for example, $X_2$ with 9 windows. As shown in Fig.1, any two points $p$ and $q$ are joined by an arc in $X_2$. Therefore, $X_2$ is arcwise connected, that is, $X_2$ is connected[8]. In the same way, each $X_i$ is convinced to be connected and compact. Applying 2-F. to a compact normal space $[0,1]\times [0,1]$, we can see that non-empty compact set $\displaystyle\bigcap _{i\in \bf{N}}X_i= $SC is connected. There exists a fundamental difference in topological nature between SC and CMTS. CMTS is 0-dim, i.e, totally disconnected, on the other hand, SC is connected. The procedure of construction of the self-similar quotient space of CMTS is based on the existence of an adequate partition of CMTS and the existence of this partition is based on the 0-dim of CMTS. Thus, concerning SC which is connected, no non-trivial quotient space which is self-similar is obtained in an elegant procedure which is employed for CMTS.\\

\bigskip

\begin{figure}[!h]
\begin{center}
\includegraphics[width=5cm,clip]{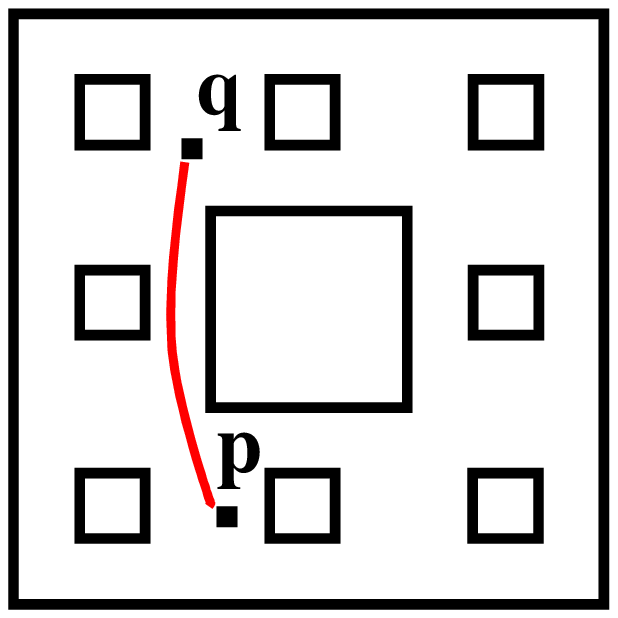}
\caption{The second step $X_2$ of the construction of Sierpi\'nski carpet$=\displaystyle\bigcap _{i\in \bf{N}}X_i$. Any two points $p$ and $q$ are joined by an arc in $X_2$.}
\end{center}

\label{f1}
\end{figure}

\section{Conclusion}

Results obtained concerning the difference in topological nature between CMTS and SC are as follows.\\

\noindent
1)~~A topological space which is homeomorphic to a self-similar space is self-similar. There exists a non-trivial quotient space $\mathcal{D}$ of CMTS which is homeomorphic to CMTS. Since CMTS is self-similar, quotient space $\mathcal{D}$ is also self-similar.\\

\noindent
2)~~Above result 1) owes mainly to the 0-dim (the totally disconnectedness) of the structures of CMTS. Therefore, it is impossible for SC which is a connected space to have a self-similar quotient space under the exactly same elegant procedure of construction as that employed for CMTS.\\ 

\noindent
~~Finally we note that, because of the connectedness of their structures, the same situations with Sierpi\'nski carpet exist for Sierpi\'nski gasket and Sierpi\'nski sponge (Menger sponge)[9].\\
~~\\

\section*{Appendix}

\noindent
1).~~Proof of the statement in 2-B.\\
~~Let a topological space $(X, T)$ be self-similar and homeomorphic to $(Y, \tau)$. By the above definition of self-similarity, there exists a metric $d$ on $X$ such that $T = \tau_d$ and there exists a system of contractions \\ $p_j:(X,\tau_d)\to(X,\tau_d)$, $d(p_j(x),p_j(x'))\le\alpha_j(\eta)d(x,x')$ for $d(x,x')<\eta,~0\le\alpha_j(\eta)<1,~j=1,\ldots,m~(2\le m<\infty)$ satisfying the relation $\displaystyle\bigcup_{j=1}^m p_j(X)=X$. Using a homeomorphism $h:(X,\tau_d)\rightarrow (Y,\tau)$, we can define a metric $\rho$ on $Y$ as 
$$
\rho(y,y')=d(h^{-1}(y),h^{-1}(y')),~~y,~y'\in Y.
$$
The metric topology $\tau_\rho$ is identical with the initial topology $\tau$. 
From the relations 1) and 2) below, $(Y, \tau) (= (Y,\tau_\rho))$ is self-similar based on a system of contractions  $q_j:(Y,\tau_\rho)\to(Y,\tau_\rho),~j=1,\ldots,m$ where $q_j$ is topologically conjugate to $p_j$ with the above homeomorphism $h$, that is, $q_j=h\circ p_j\circ h^{-1}$ .
\begin{itemize}
\item[1)]$\rho(q_j(y),q_j(y'))=d(h^{-1}(q_j(y)),h^{-1}(q_j(y')))$\\
$=d(p_j(h^{-1}(y)),p_j(h^{-1}(y')))\le\alpha_j(\eta)d(h^{-1}(y),h^{-1}(y'))$\\
$=\alpha_j(\eta)\rho(y,y')$ for $\rho(y,y')<\eta$.
\item[2)]$\displaystyle\bigcup_{j=1}^mq_j(Y)=\bigcup_{j=1}^mq_j(h(X))=h(\bigcup_{j=1}^mp_j(X))=h(X)=Y.$  ~$\Box$
\end{itemize}
~~\\

\noindent
2).~~Proof of the statement in 2-C.\\
~~Let $C(Z)$ and $d_H$ be the set of all non-empty compact sets of $(Z,\tau_d)$ and the Hausdorff metric, respectively. It is known[10] that the set dynamical system
\begin{center}
$T : (C(Z),\tau_{d_H})\rightarrow (C(Z),\tau_{d_H}),~A\mapsto \displaystyle\bigcup_{j=1}^m f_j(A)$
\end{center}
has unique fixed point $X$; that is, there exists unique non-empty compact set $X(\subset Z)$ such that $\displaystyle\bigcup_{j=1}^m f_j(X)=X.$\\
~~The condition $iii)$ is known[11] as a condition which makes the set $X$ 0-dim.\\~~To complete the proof, let us show that the conditions $i)$ and $ii)$ are sufficient for the subspace $(X,\tau_{d_X})$ of $(Z,\tau_d)$ ($d_X$ denotes the restriction of the metric $d$ on $X$) to be perfect. Namely, let us show that any point $x$ of $X$ is an accumulation point of $X$. Since any contraction on a complete metric space has one, only one, fixed point[10], it is obvious that the condition $ii)$ requires that $X$ is not a singleton. Let $x_0$ and $x_0'$ be different two points of $X$. For any point $x\in X$ and for any $\varepsilon >0$, let us consider an open sphere $S(x,\varepsilon )=\{z\in Z;d(x,z)<\varepsilon \}$. The relation $\displaystyle\bigcup_{j{_1}\cdot \cdot \cdot j{_k}}X_{j{_1}\cdot \cdot \cdot j{_k}}=X$ is valid for any $k\geq 1$[11]. Here $X_{j{_1}\cdot \cdot \cdot j{_k}}$ denotes $f_{j_1}\circ \cdot \cdot \cdot \circ f_{j_k}(X),~j_1\in \{1,\dots,m\},\dots,j_k\in \{1,\dots,m\}$, and the diameter of $X_{j{_1}\cdot \cdot \cdot j{_k}} \rightarrow 0~(k\rightarrow \infty )$[11]. There exist $k$ and a $k$-tuple $j_1\cdot \cdot \cdot j_k$ such that $x\in X_{j_1\cdot \cdot \cdot j_k}$ and the diameter of $X_{j_1\cdot \cdot \cdot j_k}<\varepsilon $. Since the contractions $f_1,\dots,f_m$ are all one to one from the condition $i)$, the point $p=f_{j_1}\circ \cdot \cdot \cdot \circ f_{j_k}(x_0)\in X_{j_1\cdot \cdot \cdot j_k}\subset X$ and the point $p'=f_{j_1}\circ \cdot \cdot \cdot \circ f_{j_k}(x_0')\in X_{j_1\cdot \cdot \cdot j_k}\subset X$ must be different. Therefore, at least either $S(x,\varepsilon )-\{x\}\ni p$ or $S(x,\varepsilon )-\{x\}\ni p'$ holds. This means that the point $x$ is an accumulation point of $X$.~$\Box$\\
~~\\

\noindent
3).~~Proof of the statement in 2-E.\\
~~Since $(Z,\tau)$ is perfect, $(Z,\tau)$ has at least two points $z$ and $z'$. Without loss of generality, there exists $U\in\tau$ such that $z\in U$ and $z'\not\in U$. 0-dim guarantees the existence of $u\in\tau\cap \Im$  such that $z\in u\subset U$. Taking this open and closed set $u$ as $Y$, we obtain a 0-dim, perfect, proper subspace $(Y,\tau_Y)$ of $(Z,\tau)$.~$\Box$\\
~~\\

\section*{Acknowledgement}
\noindent
\rm ~The authors are grateful to Prof. Emeritus H. Fukaishi, Dr. H. Ryo, Dr. Y. Yamashita, Dr. K. Yoshida, Dr. Y. Watayoh and Dr. K. Shikama for helpful discussions.   
~~\\
~~\\
\vspace{-4mm}

\section*{References}
\vspace{-6mm}

~~$\\$
[1] S.K.Ma, Modern~theory~of~critical~phenomena, Benjamin (1976) p.246.$\\$
[2] A.Fern\'andez, J. Phy. A {\bf 21}, L295 (1988).$\\$
[3] A.Illanes and S.B.Nadler Jr., Hyperspace, Marcel Dekker (1999) p.102.$\\$  
[4] A.Kitada, Isohkukan to sono ouyo, Asakura shoten (2011) p.39 (Japanese).$\\$ 
[5] The map $h:(A,\tau)\rightarrow (A',\tau')$ is said to be a homeomorphism provided that $h$ is one to one, onto and continuous together with the inverse function $h^{-1}:(A',\tau')\rightarrow (A,\tau)$. Roughly speaking, "the space $A$ is homeomorphic to the space $A'$" means that $A$ and $A'$ are geometrically equivalent. Ref.[4] p.24.$\\$
[6] S.B.Nadler Jr.,  Continuum theory, Marcel Dekker (1992) p.105, p.109.$\\$
[7] S.B.Nadler Jr.,  Continuum theory, Marcel Dekker (1992) p.6.$\\$
[8] A.Kitada, Isohkukan to sono ouyo, Asakura shoten (2011) p.48  (Japanese).$\\$
[9] M.Matsushita, Furakutal no butsuri(I), shoukabou (2004) p.15 (Japanese).$\\$
[10] S.Nakamura, T.Konishi and A.Kitada, Journal of the Physical Society of Japan {\bf 64} (1995) 731.$\\$
[11] A.Kitada, T.Konishi, T.Watanabe, Chaos, Solitons $\&$ Fractals {\bf 13} (2002) 363.$\\$

\end{document}